\documentstyle[12pt]{article}
\newcommand{\qed}{\hfill\rule{4pt}{8pt}\par\vspace{\baselineskip}}

\newtheorem{de}{Definition}[section]

\newtheorem{pr}[de]{Proposition}
\newtheorem{co}[de]{Corollary}
\newtheorem{re}[de]{Remark}

\newtheorem{te}[de]{Theorem}

\def\Box{\mbox{$\sqcap\!\!\!\!\sqcup$}}

\def\ot{\otimes}
\def\ra{\rightarrow}

\def\al{\alpha}
\def\veps{\varepsilon}

\def\bea{\begin{eqnarray*}}
\def\eea{\end{eqnarray*}}

\begin{document}
\title{On the dimension of the space of integrals on coalgebras}
\author{
S. D\u{a}sc\u{a}lescu, C. N\u{a}st\u{a}sescu and B. Toader \\[2mm]
University of Bucharest, Facultatea de Matematica,\\ Str.
Academiei 14, Bucharest 1, RO-010014, Romania,\\
e-mail: sdascal@fmi.unibuc.ro, Constantin\_nastasescu@yahoo.com,\\
bogdantd@yahoo.com }

\date{}
\maketitle

\begin{abstract}
We study the injective envelopes of the simple right
$C$-comodules, and their duals, where $C$ is a coalgebra. This is
used to give a short proof and to extend a result of Iovanov on
the dimension of the space of integrals on coalgebras. We show
that if $C$ is right co-Frobenius, then the dimension of the space
of left $M$-integrals on $C$ is $\leq {\rm dim}M$ for any left
$C$-comodule $M$ of finite support, and the dimension of the space
of right $N$-integrals on $C$ is $\geq {\rm dim}N$ for any right
$C$-comodule $N$ of finite support. If $C$ is a coalgebra, it is
discussed how far is the dual algebra $C^*$ from being
semiperfect. Some examples of integrals are computed for
incidence coalgebras.
\end{abstract}

\section{Introduction and preliminaries}
Integrals have played a key role in the structure and
representation theory of Hopf algebras. The definition of
integrals on Hopf algebras was given by Larson and Sweedler in
\cite{ls}. Answering a question posed by Sweedler, Sullivan proved
the uniqueness of integrals, i.e. that the dimension of the space
of left (or right) integrals on a Hopf algebra is either 0 or 1.
It turned out that the existence of integrals is closely related
to some coalgebraic properties of the Hopf algebra. We recall that
a coalgebra $C$ is called right semiperfect if the category ${\cal
M}^C$ of right $C$-comodules has enough projectives, or
equivalently the injective envelope of any simple left
$C$-comodule is finite dimensional. $C$ is called right
co-Frobenius if $C$ embeds in $C^*$ as a right $C^*$-module. We
have similar notions to the left. The semiperfect property and the
co-Frobenius property are not left-right symmetric. Also, if $C$
is right (left) co-Frobenius, then $C$ is right (left)
semiperfect, while the converse is not true in general. However,
for a Hopf algebra $H$, we have that $H$ is right semiperfect
$\Leftrightarrow$ $H$ is right co-Frobenius $\Leftrightarrow$ $H$
is left semiperfect $\Leftrightarrow$ $H$ is left co-Frobenius.
Moreover, these are equivalent to $H$ having non-zero left (or
right) integrals (see \cite{lin}).

We recall that a left integral on a Hopf algebra $H$ (over a field
$k$) is an element $T$ in the dual space $H^*$ such that
$h^*T=h^*(1)T$ for any $h^*\in H^*$. It is a simple, but very
useful remark of Doi \cite{doi} that such a $T$ is in fact a
morphism of right $H^*$-modules from $H$ to $k$. This suggests
that one may consider integrals in a more general framework: if
$C$ is a coalgebra, then for any left $C$-comodule $M$, a left
$M$-integral on $C$ is a morphism of left $C$-comodules (or
equivalently of right $C^*$-modules) from $C$ to $M$, and for any
right $C$-comodule $N$, a right $N$-integral on $C$ is a morphism
of right $C$-comodules (or equivalently of left $C^*$-modules)
from $C$ to $N$. We denote by $\int_{l,C,M}$ (respectively
$\int_{r,C,N}$) the space of left $M$-integrals (respectively
right $N$-integrals) on $C$. Using this point of view and a
homological approach, \c{S}tefan \cite{st} proved that ${\rm
dim}\int_{l,C,M}={\rm dim}\; M$ for any finite dimensional left
comodule over a right co-Frobenius coalgebra $C$ which is either
finite dimensional, or cosemisimple, or the underlying coalgebra
structure of a Hopf algebra. In particular, \c{S}tefan gave a
short proof for the uniqueness of integrals in a Hopf algebra.
Note that in all cases mentioned in the result above, $C$ is left
and right co-Frobenius. In \cite{dnt} it is proved that for any
right co-Frobenius coalgebra $C$ which is also left semiperfect,
we have that ${\rm dim}\int_{l,C,M}\leq{\rm dim}\; M$ for any
finite dimensional left $C$-comodule $M$. In \cite{i2} Iovanov
succeeded to prove this result in the case where $C$ is only right
co-Frobenius (thus to drop the left semiperfect condition), and
moreover to prove that for such a $C$ we also have that ${\rm
dim}\;N\leq {\rm dim}\int _{r,C,N}$ for any finite dimensional
right $C$-comodule $N$.

The main aim of this paper is to present a new approach to the
dimension of the space of left or right integrals on a right
co-Frobenius coalgebra, which produces a short proof of Iovanov's
result, and in fact slightly extends this result. If $C$ is a
coalgebra, we say that a left $C$-comodule $M$, with comodule
structure map $\rho:M\ra C\ot M$, has finite support if there
exists a finite dimensional subspace $X$ of $C$ such that $\rho
(M)\subseteq X\ot M$ (this is equivalent to the fact that the
coalgebra of coefficients of $M$ is finite dimensional). We prove
in Section \ref{seccoF}  that if $C$ is right co-Frobenius, then
${\rm dim}\int _{l,C,M}\leq {\rm dim}\;M$ (as cardinal numbers)
for any left $C$-comodule $M$ of finite support, and ${\rm
dim}\;N\leq {\rm dim}\int _{r,C,N}$ for any right $C$-comodule $N$
of finite support. Our approach to this result has as a main theme
the study of the subspaces ${\rm Rat}(_{C^*}C^*)$, ${\rm
Rat}(C^*_{C^*})$, $\oplus_{j\in J}E(T_j)^*$ and $\oplus_{i\in
I}E(S_i)^*$ of $C^*$, as well as the relations between these
subspaces, where $C_0=\oplus_{i\in I}S_i$ (respectively
$C_0=\oplus_{j\in J}T_j$) are decompositions of the coradical
$C_0$ of $C$ as a direct sum of left (respectively right)
$C$-subcomodules, and by $E(M)$ we denote the injective envelope
of a (left or right) $C$-comodule in the category of comodules. In
the case where $C$ is left and right semiperfect, these four
spaces are equal. This study is done in Section \ref{secinjind}.
In the same section we explain how far is the dual $C^*$ of a
coalgebra from being a semiperfect ring. By using the rich
topological structure of $C^*$ (more precisely $C^*$ is a complete
topological ring with a certain topology determined by the
coradical filtration of $C$), we show that idempotents lift modulo
the Jacobson radical of $C^*$ for any $C$. As a consequence, $C^*$
is a semiperfect ring if and only the coradical of $C$ is finite
dimensional. As a byproduct of our approach, we obtain in Section
\ref{seccoF} that a right quasi-co-Frobenius coalgebra which is
also a basic coalgebra is necessarily right co-Frobenius, and as
an application we show that any right quasi-co-Frobenius coalgebra
is Morita-Takeuchi equivalent to a right co-Frobenius coalgebra.
In Section \ref{secinc} we give some examples by computing the
dimension of the spaces of left and right integrals on incidence
coalgebras for certain comodules. It is known that incidence
coalgebras are rarely right (or left) co-Frobenius, more precisely
only when the underlying ordering relation is equality. We show
that for a finite dimensional non-zero right comodule $M$ over an
incidence coalgebra $C$, the space of right $M$-integrals on $C$
may be infinite dimensional, and also for any pair $m,n$ of
positive integers we can find examples such that $M$ has dimension
$m$, and the space of right $M$-integrals on $C$ has dimension
$n$.

We work over a fixed field $k$.  If $C$ is a coalgebra, then the
comodule structure map of a left $C$-comodule $M$ is denoted by
$\rho:M\ra C\ot M$, $\rho (m)=\sum m_{-1}\ot m_0$. Such a $M$ is a
right $C^*$-module with action of $c^*\in C^*$ on $m\in M$ denoted
by $m\cdot c^*$. If $R$ is a ring, and $M,N$ are right
(respectively left) $R$-modules, the set of morphisms of right
(respectively left) $R$-modules from $M$ to $N$ is denoted by
${\rm Hom}_{-R}(M,N)$ (respectively ${\rm Hom}_{R-}(M,N)$). This
notation will distinguish which side we are working in the case
where $M$ and $N$ are $R$-bimodules. For basic definitions and
facts about coalgebras we refer to \cite{dnr}, while for general
facts about modules to \cite{af}.

\section{Injective indecomposable comodules and their duals}
\label{secinjind}

Let $C$ be a coalgebra. Then the coradical $C_0$ of $C$ is the
socle of the right $C$-comodule $C$, and also the socle of the
left $C$-comodule $C$ (see \cite[Proposition 3.1.4]{dnr}). We
write $C_0=\oplus_{j\in J}T_j$, as a direct sum of simple right
$C$-comodules. Then for any $j\in J$ there exists an injective
envelope $E(T_j)$ of $T_j$ (in the category ${\cal M}^C$) such
that $E(T_j)$ is a right subcomodule of $C$ and $C=\oplus_{j\in
J}E(T_j)$ (see \cite{green} or \cite[Theorem 2.4.16]{dnr}). We
have that $E(T_j)$ is an indecomposable injective object in ${\cal
M}^C$. Then $C^*\simeq \prod_{j\in J}E(T_j)^*$ as right
$C^*$-modules. In fact we identify $C^*$ with $\prod_{j\in
J}E(T_j)^*$, by regarding $E(T_j)^*$ as the set of all elements
$c^*\in C^*$ such that $c^*(E(T_p))=0$ for any $p\in J-\{j\}$.
Then we can also consider the right $C^*$-submodule $\oplus_{j\in
J}E(T_j)^*$ of $C^*$, which is a dense subspace of $C^*$ in the
finite topology (see for example \cite[Exercise 1.2.17]{dnr}).

Let us note that the subspace $\oplus_{j\in J}E(T_j)^*$ of $C^*$
does not depend on the representation $C_0=\oplus_{j\in J}T_j$ as
a direct sum of simple right comodules. Indeed, for each
isomorphism type of simple right $C$-comodule there are only
finitely many $T_j$'s of that type, say $T_{j_1},\ldots,T_{j_n}$,
and $T_{j_1}\oplus\ldots\oplus T_{j_n}$ is just the (simple)
subcoalgebra $D$ of coefficients associated to that type of simple
comodule. Then $E(T_{j_1})^*\oplus \ldots \oplus
E(T_{j_n})^*=E(D)^*$, where $E(D)$ is the injective envelope of
the right $C$-comodule $D$. Thus $\oplus_{j\in
J}E(T_j)^*=\oplus_{\lambda}E(D_{\lambda})^*$, where
$C_0=\oplus_{\lambda}D_{\lambda}$ is the decomposition of $C_0$ as
a direct sum of simple subcoalgebras (and the $D_{\lambda}$'s are
uniquely determined by $C_0$).

For any $r\in J$ let $\veps_r\in E(T_r)^*\subseteq C^*$ such that
${\veps_r}_{|E(T_r)}=\veps$, and ${\veps_r}_{|E(T_j)}=0$ for any
$j\neq r$. It is easy to see that $\veps_rc^*=c^*$ for any $c^*\in
E(T_r)^*$, in particular $\veps_r^2=\veps_r$, and also that
$\veps_rc^*=0$ for any $c^*\in E(T_j)^*$, with $j\neq r$. As a
consequence, we have that $E(T_r)^*=\veps_rC^*$.

The following will be a key result in the sequel.

\begin{pr} \label{isofinsup}
Let $M$ be a left $C$-comodule. Then the map
$$\gamma_M:M\ra Hom_{-C^*}(\oplus_{j\in
J}E(T_j)^*,M),\; \gamma_M(m)(c^*)=m\cdot c^*$$ is an injective
linear map. Moreover, if $M$ has finite support, then $\gamma_M$
is a linear isomorphism.
\end{pr}
{\bf Proof:} Let $m\in M$ such $\gamma_M(m)=0$. Write
$\rho(m)=\sum m_{-1}\ot m_0$, where $\rho$ is the comodule
structure map of $M$. Since $\oplus_{j\in J}E(T_j)^*$ is dense in
$C^*$, there exists $c^*\in \oplus_{j\in J}E(T_j)^*$ which is
equal to $\veps$ on all $m_{-1}$'s. Then $m=m\cdot \veps=m\cdot
c^*=\gamma_M(m)(c^*)=0$, so $\gamma_M$ is injective.

Assume now that $M$ has finite support. Then there exists a finite
subset $J_0$ of $J$ such that $\rho(M)\subseteq (\oplus_{j\in
J_0}E(T_j))\ot M$. Then for any $m\in M$ we have that $m\cdot
(\sum_{j\in J_0}\veps_j)=m$ and $m\cdot c^*=0$ for any $c^*\in
E(T_r)^*$ with $r\notin J_0$. Let $\phi \in
Hom_{-C^*}(\oplus_{j\in J}E(T_j)^*,M)$, and let $m=\sum_{j\in
J_0}\phi(\veps_j)$. We show that $\phi =\gamma_M(m)$, which will
prove that $\gamma_M$ is also
surjective. \\
If $c^*\in E(T_r)^*$ with $r\notin J_0$, then
$\phi(\veps_r)=\phi(\veps_r^2)=\phi(\veps_r)\cdot\veps_r=0$, hence
\bea
\phi(c^*)&=&\phi(\veps_rc^*)\\&=&\phi(\veps_r)c^*\\&=&0\\&=&m\cdot
c^*\\&=&\gamma_M(m)(c^*)\eea If $c^*\in E(T_r)^*$ with $r\in J_0$,
then \bea
\gamma_M(m)(c^*)&=&m\cdot c^*\\
&=&\sum_{j\in J_0}\phi(\veps_j)\cdot c^*\\
&=&\sum_{j\in J_0}\phi(\veps_jc^*)\\
&=&\phi(c^*)\eea the last equality holding since $\veps_jc^*=0$
for any $j\neq r$, and $\veps_rc^*=c^*$.\\
We conclude that $\phi =\gamma_M(m)$, which ends the proof. \qed

\begin{pr} \label{inclratst}
With notation as above we have that ${\rm Rat}(_{C^*}C^*)\subseteq
\oplus_{j\in J}E(T_j)^*$.
\end{pr}
{\bf Proof:} Let $c^*\in {\rm Rat}(_{C^*}C^*)$. Then there exist
finite families $(c_{\al})_{\al}\subseteq C$ and
$(c^*_{\al})_{\al}\subseteq C^*$ such that
$d^*c^*=\sum_{\al}d^*(c_{\al})c^*_{\al}$ for any $d^*\in C^*$.
Choose a finite subset $F$ of $J$ such that all $c_{\al}$'s lie in
$\oplus_{j\in F}E(T_j)$.

Let $d^*\in C^*$ such that $d^*_{|E(T_r)}=\veps$ for any $r\notin
F$, and $d^*_{|E(T_r)}=0$ for any $r\in F$. Then $d^*(c_{\al})=0$
for any $\al$, so then $d^*c^*=\sum_{\al}d^*(c_{\al})c^*_{\al}=0$.
On the other hand, if $r\notin F$ and $c\in E(T_r)$, then
$\Delta(c)\in E(T_r)\ot C$, so then $(d^*c^*)(c)=\sum
d^*(c_1)c^*(c_2)=\sum \veps (c_1)c^*(c_2)=c^*(c)$. Since
$d^*c^*=0$ we get that $c^*(c)=0$, and this shows that $c^*\in
\oplus_{r\in F}E(T_r)^*\subseteq \oplus_{j\in J}E(T_j)^*$. \qed

Let us note that ${\rm Rat}(_{C^*}C^*)$ is a submodule of the left
$C^*$-module $C^*$, while $\oplus_{j\in J}E(T_j)^*$ is a submodule
of the right $C^*$-submodule $C^*$.

Similarly, if we work with left $C$-comodules, we can write
$C_0=\oplus_{i\in I}S_i$, a direct sum of simple left comodules,
then $C=\oplus_{i\in I}E(S_i)$, a direct sum of indecomposable
injective left $C$-comodules, and then we have that ${\rm
Rat}(C^*_{C^*})\subseteq \oplus_{i\in I}E(S_i)^*$, where ${\rm
Rat}(C^*_{C^*})$ is a submodule of the right $C^*$-module $C^*$,
while $\oplus_{i\in I}E(S_i)^*$ is a submodule of the left
$C^*$-submodule $C^*$, which does not depend on the choice of the
simples $(S_i)_i$ in the decomposition $C_0=\oplus_i S_i$. The
following result shows that under certain finiteness conditions on
$C$, the subspaces ${\rm Rat}(_{C^*}C^*)$, ${\rm Rat}(C^*_{C^*})$,
$\oplus_{j\in J}E(T_j)^*$ and $\oplus_{i\in I}E(S_i)^*$ are in a
special relation.

\begin{pr} \label{inclrightsemiperfect}
Let $C$ be a right semiperfect coalgebra. Then
$${\rm Rat}(C^*_{C^*})\subseteq \oplus_{i\in I}E(S_i)^*\subseteq {\rm Rat}(_{C^*}C^*)\subseteq
\oplus_{j\in J}E(T_j)^*$$
\end{pr}
{\bf Proof:} Since $C$ is right semiperfect, each $E(S_i)$ is
finite dimensional. Then $E(S_i)^*$ is a rational left
$C^*$-module (see for example \cite[Lemma 2.2.12]{dnr}), so
$E(S_i)^*\subseteq {\rm Rat}(_{C^*}C^*)$. Now all follows from
Proposition \ref{inclratst} and its left hand side version. \qed

\begin{re} There is, of course, a version of Proposition
\ref{inclrightsemiperfect} for left semiperfect coalgebras. These
two results show that if $C$ is left and right semiperfect, then
$${\rm Rat}(_{C^*}C^*)={\rm Rat}(C^*_{C^*})=\oplus_{i\in
I}E(S_i)^*=\oplus_{j\in J}E(T_j)^*$$ This fact is already known
(see \cite[Theorem 2.4]{bdgn}).\end{re}

Part of the following result appears  in \cite[Lemma 1.4]{i1}. For
completeness we include a proof, which is new and seems to be
shorter than the one in \cite{i1}.

\begin{pr} \label{dualindec}
Let $C$ be a coalgebra, $T$  a simple right $C$-comodule and
$E(T)$ its injective envelope in the category ${\cal M}^C$. Then
$E(T)^*$ is an indecomposable right $C^*$-module and
$End_{-C^*}(E(T)^*)$ is a local ring. Moreover $E(T)^*$ is a
projective cover of the simple right $C^*$-module $T^*$.
\end{pr}
{\bf Proof:} It is enough to prove the statement for $T=T_j$,
where $j\in J$. Let $j\in J$. The dual of the inclusion morphism
$T_j\ra E(T_j)$ is a surjective morphism of right $C^*$-modules
whose kernel is $T_j^{\perp}=\{ g\in E(T_j)^*|\; g(T_j)=0\}$. Then
we have an isomorphism of right $C^*$-modules
$E(T_j)^*/T_j^{\perp}\simeq T_j^*$. Since clearly $T_j^*$ is a
simple right $C^*$-module, we obtain that $T_j^{\perp}$ is a
maximal right submodule of $E(T_j)^*$. Therefore ${\rm
Rad}(E(T_j)^*)\subseteq T_j^{\perp}$, where by ${\rm Rad}(M)$ we
denote the Jacobson radical of a right $C^*$-module $M$.

It is known (see for example \cite[Proposition 3.1.8]{dnr}) that
$${\rm Rad}(C^*)=C_0^{\perp}=\{c^*\in C^*|\; c^*(C_0)=0\}$$ Since
$C_0=\oplus _{j\in J}T_j$ and $C^*$ is identified with
$\prod_{j\in J}E(T_j)^*$, we have that
$$\prod_{j\in J}T_j^{\perp}=C_0^{\perp}={\rm Rad}(C^*)={\rm Rad}(\prod_{j\in J}E(T_j)^*)
\subseteq \prod_{j\in J}{\rm Rad}(E(T_j)^*)\subseteq \prod_{j\in
J}T_j^{\perp}$$ which shows that both inclusions in the sequence
of relations in the row above are equalities. In particular ${\rm
Rad}(E(T_j)^*)= T_j^{\perp}$, thus $T_j^{\perp}$ is the unique
maximal submodule of $E(T_j)^*$, and then clearly $T_j^{\perp}$ is
superfluous in $E(T_j)^*$. Thus $E(T_j)^*$ together the dual of
the inclusion morphism is a projective cover of $T_j^*$.

Now since $E(T_j)^*$ is projective (as a direct summand of $C^*$),
we have by \cite[Proposition 17.19]{af} that
$End_{-C^*}(E(T_j)^*)$ is a local ring, in particular $E(T_j)^*$
is indecomposable.\qed

\begin{re}
Let us note that the relation ${\rm Rad}(E(T_j)^*)= T_j^{\perp}$,
that we showed in the proof of Proposition \ref{dualindec}, is
just a local (corepresentation) version of the fact that ${\rm
Rad}(C^*)=C_0^{\perp}$.
\end{re}

Concerning the existence of projective covers in the category
${\cal M}_{C^*}$, we can prove even more than what is obtained in
Proposition \ref{dualindec}. The following result shows that any
right $C^*$-module which is the dual of a finite dimensional right
$C$-comodule has a projective cover.

\begin{pr}
Let $M$ be a finite dimensional right $C$-comodule. Then $M^*$ has
a projective cover in the category of right $C^*$-modules.
\end{pr}
{\bf Proof:} Write the socle of $M$ as $T_1\oplus \ldots \oplus
T_n$, a sum of simple right $C$-comodules. Then $E(M)=E(T_1)\oplus
\ldots \oplus E(T_n)$, so we have an epimorphism of right
$C^*$-modules $E(T_1)^*\oplus \ldots \oplus E(T_n)^*\simeq
E(M)^*\ra M^*$. Now proceed as in the proof of \cite[Theorem
27.6]{af}. Denote ${\rm Rad}(C^*)={\cal J}$, the Jacobson radical
of $C^*$. Then we have an epimorphism of right $C^*$-modules
$E(T_1)^*/E(T_1)^*{\cal J}\oplus \ldots \oplus
E(T_n)^*/E(T_n)^*{\cal J}\ra M^*/M^*{\cal J}$, which shows that
$M^*/M^*{\cal J}$ is a finite direct sum of $C^*$-modules of the
form $E(T)^*/E(T)^*{\cal J}$ (and these are all simple since by
\cite[Proposition 17.19]{af} we have that ${\rm
Rad}(E(T)^*)=E(T)^*{\cal J}=T^\perp$). By Proposition
\ref{dualindec}, each $E(T)^*/E(T)^*{\cal J}$ has a projective
cover, hence so does $M^*/M^*{\cal J}$. Let $f:P\ra M^*/M^*{\cal
J}$ be a projective cover. Then $f$ induces a morphism $g:P\ra
M^*$ such that $f=g\pi$, where $\pi:M^*\ra M^*/M^*{\cal J}$ is the
natural projection. By Nakayama's Lemma, $\pi$ is a superfluous
epimorphism, and now using \cite[Lemma 27.5]{af}, we obtain that
$g:P\ra M^*$ is a projective cover of $M^*$. \qed

We recall from \cite{af} that a ring $R$ is called semiperfect if
$R/{\rm Rad}(R)$ is semisimple and idempotents lift modulo ${\rm
Rad}(R)$ (or equivalently any simple right $R$-module has a
projective cover). We are interested to see when is the dual
algebra $C^*$ of a coalgebra $C$ a semiperfect ring. The following
shows that the second condition in the definition of
semiperfectness is always satisfied. In order to prove this, we
use the rich topological structure of $C^*$.

\begin{pr} \label{idlift}
Let $C$ be a coalgebra. Then idempotents lift modulo ${\rm
Rad}(C^*)$, the Jacobson radical of the dual algebra $C^*$.
\end{pr}
{\bf Proof:} Consider the algebra filtration $\ldots \subseteq
F_{-1}\subseteq F_0\subseteq F_1\subseteq \ldots $ of $C^*$, where
$F_{-n}=C_{n-1}^{\perp}$ for any positive integer $n$, and
$F_n=C^*$ for any non-negative integer $n$.  Then $C^*$ is a
topological ring with the family $(C_n^{\perp})_{n\geq 0}$ as a
fundamental system of neighborhoods of 0 (see for example
\cite[Chapter D, Section I]{nv}). The completion of $C^*$ in this
topology is \bea \widehat{C^*}&=&\lim_{\longleftarrow\atop
{p\in{\bf Z}}}\; C^*/F_p\\
&=&\lim_{\longleftarrow\atop {n\geq 0}}\; C^*/C_n^{\perp} \\
&=&\lim_{\longleftarrow\atop {n\geq 0}}\; C_n^*\\
&=&\lim_{\longleftarrow\atop {n\geq 0}}Hom_{{\cal M}^C}(C_n,C)\\
&=&Hom_{{\cal M}^C}(\cup_{n\geq 0}C_n,C)\\
&=&Hom_{{\cal M}^C}(C,C)\\ &\simeq&C^* \eea Therefore $C^*$ is a
complete topological ring.

Clearly any element $x\in{\rm Rad}(C^*)$ is topologically
nilpotent, i.e. the sequence $(x^n)_n$ converges to 0, since
$x^n\in (C_0^{\perp})^n\subseteq F_{-n}$ for any $n>0$. Now we can
apply \cite[Lemma VII.1, page 312]{nv} and obtain that idempotents
lift modulo ${\rm Rad}(C^*)$. \qed

Now we are able to characterize coalgebras for which the dual
algebra is a semiperfect ring.

\begin{pr}
Let $C$ be a coalgebra. Then the following assertions are
equivalent.\\
{\rm (i)} The dual algebra $C^*$ is a semiperfect ring.\\
{\rm (ii)} The coradical $C_0$ of $C$ is finite dimensional.\\
{\rm (iii)} There exist finitely many isomorphism types of simple
right (or left) $C$-comodules.
\end{pr}
{\bf Proof:} Write $C_0=\oplus_{\lambda\in \Lambda}D_{\lambda}$, a
direct sum of simple subcoalgebras. Then $C_0^*\simeq
\prod_{\lambda\in \Lambda}D_{\lambda}^*$, a direct product of
simple algebras. Therefore $C_0^*$ is a semisimple algebra if and
only if $\Lambda$ is finite. Then (i)$\Rightarrow$(ii) follows,
since if $C^*$ is semiperfect, we have that $C^*/{\rm
Rad}(C^*)=C^*/C_0^{\perp}\simeq C_0^*$ is semisimple, so $\Lambda$
must be finite, and then $C_0$ is finite dimensional. Also,
(ii)$\Rightarrow$(i) holds, since $C^*/{\rm Rad}(C^*)\simeq C_0^*$
is semisimple, and by Proposition \ref{idlift} idempotents lift
modulo the Jacobson radical of $C^*$.

The equivalence of (ii) and (iii) is clear, since each isomorphism
type of a simple right $C$-comodule is associated to a unique
simple subcoalgebra of $C$ (its coalgebra of coordinates). \qed

We have recalled that a ring is semiperfect if and only if any
simple right module over that ring has a projective cover. If $C$
is an arbitrary coalgebra, we can describe all simple right
$C^*$-modules having a projective cover.

\begin{pr}
Let $C$ be a coalgebra. Then a simple right $C^*$-module $M$ has a
projective cover if and only if $M$ is a rational $C^*$-module
(i.e. the module structure comes from a simple left $C$-comodule
structure on $M$).
\end{pr}
{\bf Proof:} Assume that $M$ has a projective cover $P$ as a right
$C^*$-module. Then by \cite[Proposition 17.19]{af} we have that
$P{\rm Rad}(C^*)$ is the unique maximal submodule of $P$, thus
${\rm Rad}(P)=P{\rm Rad}(C^*)$, $P$ is cyclic and $M\simeq P/{\rm
Rad}(P)$. Then $P\oplus Q\simeq C^*$ as right $C^*$-modules for
some $C^*$-module $Q$. Factoring out the Jacobson radical, we
obtain that $P/{\rm Rad}(P)\oplus Q/{\rm Rad}(Q)\simeq C^*/{\rm
Rad}(C^*)$, thus $M$ embeds in $C^*/{\rm Rad}(C^*)$ as a right
$C^*$-module.

But $C^*/{\rm Rad}(C^*)\simeq C_0^*$ as right $C^*$-modules, where
$C_0^*$ is viewed as a right $C^*$-module via the algebra map
$C^*\ra C_0^*$ dual to the inclusion map $C_0\ra C$. Now write
$C_0=\oplus_{\lambda\in \Lambda}D_{\lambda}$, a direct sum of
simple subcoalgebras. Then $C_0^*\simeq \prod_{\lambda\in
\Lambda}D_{\lambda}^*$, which is even an isomorphism of right
$C^*$-modules. It is easy to see that the socle of the right
$C^*$-module $\prod_{\lambda\in \Lambda}D_{\lambda}^*$ is
$\oplus_{\lambda\in \Lambda}D_{\lambda}^*$, which is a rational
right $C^*$-module, since each $D_{\lambda}^*$ is rational, as the
dual of a finite dimensional comodule (see \cite[Lemma
2.2.12]{dnr}). Therefore $M$ is rational, too.

For the converse, assume that $M$ is a simple right $C^*$-module
which is rational, so $M$ is a simple left $C$-comodule. Then
$S=M^*$ is a rational simple left $C^*$-module, and $M\simeq S^*$
as right $C^*$-modules. Now everything follows from the fact that
$S^*$ has the projective envelope $E(S)^*$ by Proposition
\ref{dualindec}. \qed

\section{Finiteness conditions and integrals on coalgebras}
\label{seccoF}

In this section we investigate the connection between the
injective indecomposable left comodules and the duals of the
injective indecomposable right comodules, and apply it to the
study of integrals on co-Frobenius coalgebras.

\begin{pr} \label{phi}
Assume that $C$ is a right co-Frobenius coalgebra. Then there
exists an injective map $\phi:I\ra J$ such that $E(S_i)\simeq
E(T_{\phi (i)})^*$ as right $C^*$-modules for any $i\in I$.
\end{pr}
{\bf Proof:} Since $C$ is right co-Frobenius, it embeds in $C^*$
as a right $C^*$-module. In fact this is an embedding of $C$ in
${\rm Rat}(C^*_{C^*})$. Since $C$ must be right semiperfect, we
have by Proposition \ref{inclrightsemiperfect} that ${\rm
Rat}(C^*_{C^*})$ embeds in $\oplus_{j\in J}E(T_j)^*$ as a right
$C^*$-module.

Let $S_{i_1},\ldots,S_{i_n}$ be the objects of a given isomorphism
type among the $S_i$'s (they are in finite number since their sum
is the simple subcoalgebra of coefficients associated to that
type). Then $E(S_{i_1})\oplus \ldots \oplus E(S_{i_n})$ embeds as
a right $C^*$-module in $\oplus_{j\in J}E(T_j)^*$. Since all
$E(S_i)$'s are finite dimensional, $E(S_{i_1})\oplus \ldots \oplus
E(S_{i_n})$ embeds in fact in $\oplus_{j\in J_0}E(T_j)^*$ for a
finite subset $J_0$ of $J$. But each $E(S_i)$ and each $E(T_j)^*$
is a right $C^*$-module with a local endomorphism ring, so using
Azumaya's Theorem we obtain that there exist $j_1,\ldots,j_n\in
J_0$ such that $E(S_{i_1})\simeq E(T_{j_1})^*,
\dots,E(S_{i_n})\simeq E(T_{j_n})^*$. We define
$\phi(i_1)=j_1,\ldots,\phi(i_n)=j_n$. Proceeding the same for all
isomorphism types of simple left $C$-comodules, we get a map
$\phi:I\ra J$ such that $E(S_i)\simeq E(T_{\phi (i)})^*$  for any
$i\in I$. It is clear that $\phi$ is injective, since for
non-isomorphic $S_{\al}$ and $S_{\beta}$ we can not have
$\phi(\al)=\phi(\beta)$. Indeed, otherwise we would have
$E(S_{\al})\simeq E(T_{\phi(\al)})^*\simeq
E(T_{\phi(\beta)})^*\simeq E(S_{\beta})$, and then $S_{\al}$ and
$S_{\beta}$ must be isomorphic as socles of their injective
envelopes. \qed

\begin{co} \label{sumdir}
Let $C$ be a right co-Frobenius coalgebra. Then\\
{\rm (1)} $C$ is a direct summand of $\oplus_{j\in J}E(T_j)^*$ as
a
right $C^*$-module.\\
{\rm (2)} $\oplus_{i\in I}E(S_i)^*$ is a direct summand of  $C$ as
a left $C^*$-module.
\end{co}
{\bf Proof:} (1) Let $\phi:I\ra J$ be an injective map whose
existence was proved in Proposition \ref{phi}. We have that
$C=\oplus_{i\in I}E(S_i)\simeq\oplus_{i\in I}E(T_{\phi(i)})^*$,
hence
$$\oplus_{j\in J}E(T_j)^*\simeq C\oplus (\oplus_{j\in
J-\phi(I)}E(T_j)^*)$$ (2) Clearly $E(T_{\phi(i)})$ is finite
dimensional for any $i\in I$, so then $E(S_i)^*\simeq
E(T_{\phi(i)})$ as left $C^*$-modules. Therefore $$C\simeq
\oplus_{j\in J}E(T_j)=(\oplus_{i\in I}E(T_{\phi(i)}))\oplus
(\oplus_{j\in J-\phi(I)}E(T_j))\simeq (\oplus_{i\in
I}E(S_i)^*)\oplus (\oplus_{j\in J-\phi(I)}E(T_j))$$ \qed

Another interesting consequence of Proposition \ref{phi} is the
following. In the case where $C$ is a Hopf algebra, the result is
already known, see \cite[Proposition 3.7]{hai}, and also
\cite{cuadra}.

\begin{co}
Let $C$ be a right co-Frobenius coalgebra. Then for any simple
left $C$-comodule $S$, the injective envelope $E(S)$ has a unique
maximal subcomodule.
\end{co}
{\bf Proof:} Proposition \ref{phi} shows that $E(S)\simeq E(T)^*$
for a simple right $C$-comodule $T$. But $E(T)^*$ has a unique
maximal submodule by the proof of Proposition \ref{dualindec}.
\qed

Now we can prove the main result of the paper.

\begin{te}
Let $C$ be a right co-Frobenius coalgebra. Then\\
(1) For any left $C$-comodule $M$ of finite support we have that
${\rm dim}\int
_{l,C,M}\leq {\rm dim}\;M$.\\
(2) For any right $C$-comodule $N$ of finite support we have that
${\rm dim}\;N\leq {\rm dim}\int _{r,C,N}$.
\end{te}
{\bf Proof:} (1) We know from Corollary \ref{sumdir} that $C$ is a
direct summand of $\oplus_{j\in J}E(T_j)^*$ as a right
$C^*$-module. Then there is a surjective morphism of vector spaces
$$Hom_{-C^*}(\oplus_{j\in J}E(T_j)^*,M)\longrightarrow Hom_{-C^*}(C,M)=\int
_{l,C,M}$$ But $Hom_{-C^*}(\oplus_{j\in J}E(T_j)^*,M)\simeq M$ by
Proposition \ref{isofinsup}, and the result is proved.\\
(2) Similarly we have a surjective morphism of vector spaces
$$\int _{r,C,N}=Hom_{C^*-}(C,N)\longrightarrow Hom_{C^*-}(\oplus_{i\in I}E(S_i)^*,N)$$ By a
left hand side version of Proposition \ref{isofinsup}, we have
that $Hom_{C^*-}(\oplus_{i\in I}E(S_i)^*,N)\simeq N$ as vector
spaces, and the result follows. \qed

\begin{co}
Let $C$ be a left and right co-Frobenius coalgebra. Then for any
left (respectively right) $C$-comodule $M$ of finite support we
have that ${\rm dim}\int _{l,C,M}={\rm dim}\;M$ (respectively
${\rm dim}\int _{r,C,M}={\rm dim}\;M$).
\end{co}

We recall from \cite{gn} that a coalgebra $C$ is called right
quasi-co-Frobenius if $C$ embeds as a right $C^*$-module in a free
$C^*$-module. This is equivalent to $C$ being projective as a left
$C$-comodule (see \cite[Theorem 1.3]{gn}). It is known that a
right co-Frobenius coalgebra is right quasi-co-Frobenius, and a
right quasi-co-Frobenius coalgebra is right semiperfect, while the
converse assertions are not true in general. We also recall that a
coalgebra $C$ is called basic if the coradical of $C$ is a direct
sum of pairwise non-isomorphic simple left coideals (see \cite{cm}
or \cite{cg}). As a byproduct of the proof of Proposition
\ref{phi} we obtain the following.

\begin{pr} \label{propquasi}
Let $C$ be a right quasi-co-Frobenius coalgebra which is also a
basic coalgebra. Then $C$ is right co-Frobenius.
\end{pr}
{\bf Proof:} Since $C$ is basic, we have that $C_0=\oplus_{i\in
I}S_i$, a direct sum of pairwise non-isomorphic simple left
$C$-subcomodules of $C$. As $C$ is right quasi-co-Frobenius, it
embeds as a right $C^*$-module in a free $C^*$-module, and then so
does each $E(S_i)$. But $E(S_i)$ is finite dimensional, so it
embeds in fact in $(C^*)^n$ for some positive integer $n$. Hence
$E(S_i)$ embeds in ${\rm Rat}(C^*_{C^*})^n$, which itself embeds
in $\oplus_{j\in J}(E(T_j)^*)^n$. Now with the Azumaya type
argument as in the proof of Proposition \ref{phi}, we see that
$E(S_i)\simeq E(T_j)^*$ for some $j$, and in this way we obtain an
injective map $\phi$ as in Proposition \ref{phi}. Then $C=\oplus
_{i\in I}E(S_i)\simeq \oplus _{i\in I}E(T_{\phi(i)})^*$ embeds as
a right $C^*$-module in $C^*$, so $C$ is right co-Frobenius. \qed

Two coalgebras are called Morita-Takeuchi equivalent if their
categories of right (or equivalently of left) comodules are
equivalent, see \cite{tak}. It is a classical result that any
quasi-Frobenius finite dimensional algebra is Morita equivalent to
a Frobenius algebra. The following is a dual version (and a
generalization) of this result, for coalgebras of arbitrary
dimension.

\begin{co} \label{corquasi}
Any right quasi-co-Frobenius coalgebra is Morita-Takeuchi
equivalent to a right co-Frobenius coalgebra.
\end{co}
{\bf Proof:} Let $C$ be a right quasi-co-Frobenius coalgebra, and
let $D$ its basic coalgebra (see for example \cite{cg}), which is
Morita-Takeuchi equivalent to $C$. Let $F:^C{\cal M}\ra ^D{\cal
M}$ be an equivalence functor between the categories of left
comodules. Then $D$ embeds as a left $D$-comodule in a direct sum
$F(C)^{(I)}$ of copies of $F(C)$, and as $D$ is injective, it is a
direct summand in $F(C)^{(I)}$. But $F(C)$ is projective, and then
so are $F(C)^{(I)}$ and $D$, showing that $D$ is also right
quasi-co-Frobenius. Since $D$ is also basic, Proposition
\ref{propquasi} shows that it is right co-Frobenius. \qed

\begin{re}
Proposition \ref{propquasi} and Corollary \ref{corquasi} were
proved in \cite{gmn} in the particular case of coalgebras that are
left and right quasi-co-Frobenius.
\end{re}

\section{Some examples of integrals on incidence coalgebras} \label{secinc}

Let $(X,\leq )$ be a locally finite partially ordered set, i.e.
the interval $[x,y]=\{z|\; x\leq z\leq y\}$ is finite for any
$x\leq y$. Let $C$ be the incidence coalgebra of $X$, which has a
linear basis $\{ e_{x,y}|x,y\in X, x\leq y\}$, and
comultiplication $\Delta$ and counit $\varepsilon$ defined by
$$\Delta(e_{x,y})=\sum_{x\leq z\leq y}e_{x,z}\ot e_{z,y}$$
$$\varepsilon (e_{x,y})=\delta_{x,y}$$
for any $x,y\in X$ with $x\leq y$ (where $\delta_{x,y}$ means
Kronecker's delta). \\
 Denote $x^+=\{y|\;x\leq y\}$ and  $x^-=\{y|\;
y\leq x\}$ for any $x\in X$.

The only simple right (or left) $C$-subcomodules of $C$ are the
one-dimensional spaces $S_x=<e_{x,x}>$, where $x\in X$. The
injective cover of the right $C$-comodule $S_x$ is
$E_r(S_x)=<e_{x,y}|\; y\in x^+>$, while the injective cover of the
left $C$-comodule $S_x$ is $E_l(S_x)=<e_{y,x}|\; y\in x^->$. Thus
$C$ is right (respectively left) semiperfect if and only if $x^-$
(respectively $x^+$) is finite for any $x\in X$ (see for example
\cite{si}). It is known that $C$ is very rarely co-Frobenius, more
precisely it is right (or equivalently left) co-Frobenius if and
only if the order relation on $X$ is the equality (see
\cite{dnv}).

Let $M\in {\cal M}^C$ with comodule structure map $\rho$. Then
giving a morphism of right $C$-comodules $f:E_r(S_x)\ra M$ is the
same with giving a family $(m_y)_{y\in x^+}\in M$ such that

\begin{equation} \label{ecmorf}
\rho (m_y)=\sum _{z\in [x,y]}m_z\ot e_{z,y} \mbox{  for any }y\in
x^+ \end{equation} Indeed, this correspondence associates to $f$
the family $(f(e_{x,y}))_{y\in x^+}$.

Fix now some $u\in X$ and let $M=E_r(S_u)$. We will compute the
right $M$-integrals on $C$. Let $f$ be such an integral. Then for
a fixed $x\in X$, the restriction of $f$ to $E_r(S_x)$ is given by
a family $(m_y)_{y\in x^+}\in E_r(S_u)$ satisfying equation
(\ref{ecmorf}). Write $m_y=\sum_{v\in u^+}\al_{y,v}e_{u,v}$ for
some scalars $\al_{y,v}$ (only finitely many non-zero, for each
$y$). Then equation (\ref{ecmorf}) writes

\begin{equation} \label{eccalcul}
\sum _{v\in u^+}\sum_{p\in [u,v]}\al_{y,v}e_{u,p}\ot
e_{p,v}=\sum_{z\in [x,y]}\sum_{r\in u^+}\al_{z,r}e_{u,r}\ot
e_{z,y}
\end{equation}

We see that common tensor monomials in the left and right hand
sides of equation (\ref{eccalcul}) occur only for $z=p=r$ and
$v=y$. Hence if it is not true that $u\leq y$, such common terms
can not occur, and we have $\al_{y,v}=0$ for any $v\in u^+$, which
means that $m_y=0$.

Assume now that $u\leq y$. Then the only common tensor monomials
in the left and right hand sides are of the form $e_{u,p}\ot
e_{p,y}$ with coefficients $\al_{y,y}$ in the left hand side, and
$\al_{p,p}$ in the right hand side, where $p\in [u,y]\cap [x,y]$.
We obtain that $\al_{y,y}=\al_{p,p}$ for any such $p$, and
$\al_{y,v}=0$ for any $v\in u^+$, $v\neq y$. Moreover, if there
exists $p\in [u,y]-[x,y]$ (and this is clearly equivalent to the
fact that $u\notin [x,y]$), the term $e_{u,p}\ot e_{p,y}$ shows up
in the left hand side, but not in the right hand side, and we must
have $\al_{y,y}=0$, and then $m_y=0$. In the case where such a $p$
does not exist, i.e. $u\in [x,y]$, we must have
$\al_{y,y}=\al_{u,u}$, a scalar not depending on $y\in u^+$.
Therefore
$$f(e_{x,y})=\left\{
\begin{array}{l}
0,\mbox{ if }u\notin [x,y]\\
\alpha e_{u,y}, \mbox{ if } u\in [x,y]
\end{array}\right.$$
where $\alpha$ is a scalar. This shows that $$Hom_{{\cal
M}^C}(E_r(S_x),E_r(S_u))=0 \mbox{  if }x\notin u^-$$ and
$${\rm dim}\; Hom_{{\cal M}^C}(E_r(S_x),E_r(S_u))=1 \mbox{ if }x\in
u^-$$ Then
$$Hom_{{\cal M}^C}(C,E_r(S_u))\simeq \prod_xHom_{{\cal
M}^C}(E_r(S_x),E_r(S_u))=\prod_{x\in u^-}Hom_{{\cal
M}^C}(E_r(S_x),E_r(S_u))\simeq k^{u^-}$$ In conclusion, if
$M=E_r(S_u)$ we have
that \\

$\bullet$ If $u^+$ is finite and $u^-$ is infinite, then $M$ is a
finite dimensional right $C$-comodule, and $\int_{r,C,M}$ is
infinite dimensional. \\

$\bullet$ If both $u^+$ and $u^-$ are finite dimensional, then $M$
is a right $C$-comodule of dimension $|u^+|$, and ${\rm dim}\;
\int_{r,C,M}=|u^-|$.

These show that for a non-zero finite dimensional right
$C$-comodule $M$, the space $\int_{r,C,M}$ may be infinite
dimensional, and also the pair $({\rm dim} \; M,{\rm dim}
\int_{r,C,M})$ may be any pair of positive integers (for certain
choices of $X$ and $u$).

\end{document}